\documentclass[12pt]{article}

\usepackage{amssymb,amsmath}

\usepackage[francais,english]{babel}
\usepackage[T1]{fontenc}

 \makeatletter \@addtoreset{equation}{section}

\usepackage{hhline} 

\usepackage{graphicx}

 \makeatletter \@addtoreset{equation}{section}
\newtheorem{thm}{Theorem}[section]
\newtheorem{lem}{Lemma}[section]
\newtheorem{cor}{Corollary}[section]
\newtheorem{prop}{Proposition}[section]
\newtheorem{rem}{Remark}[section]
\newtheorem{defi}{Definition}[section]

\begin{document}


\title{Computation of maximal turning points \\ by a variational approach}


\author{Yavdat Il'yasov  \\Institute of Mathematics RAS,
Ufa, Russia \and Alexsandr Ivanov \\ Bashkir State University, Ufa, Russia}
\date{$\ $}

\maketitle

\begin{abstract}
We develop a new approach for finding bifurcations of solutions of nonlinear problems, which is based on the detection of extreme values of a new type of variational functional associated with the considered problems. The variational functional  is obtained constructively by the extended functional method which can be applied to a wide class of parametric problems including nonlinear partial differential equations.    Sufficient and necessary conditions for the existence of a maximal turning point by the approach are proved.
Based on these, an algorithm  of the quasi-direction of steepest ascent is introduced. 
Simulation experiments are used to illustrate the behaviour of the method and to discuss  its  advantages and disadvantages in comparison with the alternatives.
\end{abstract}


\section{Introduction}

The present paper is devoted to determining and computation of turning point type bifurcations of solutions branches  $u_\lambda \in \mathbb{R}^n$, $ \lambda \in (a,b)$ to systems of nonlinear equations of the form
\begin{equation}
\tag{$\mathcal{F}$}
\label{F}
F(u,\lambda):=T(u)-\lambda G(u)=0,
\end{equation}
where $T,G: \mathbb{R}^n \to \mathbb{R}^n$ are continuously differentiable functions of $u\in \mathbb{R}^n$. 
The generation of the branch of solutions and drawing the associated bifurcation diagram
is important in the investigation of various mathematical models arising in  physics, control theory, biology, ecology, economics and many other areas of science and technology and has always attracted
attention of researchers see e.g. \cite{Doedel,keller1,keller2,seydel}.

In the modern literature on the numerical analysis, to compute the solution curves and their turning points are commonly used continuation methods (see e.g. \cite{doedK, keller1,  seydel}). In general, the continuation algorithm consists of (i) selecting the starting point $(u_0,\lambda_0)$ which belongs to (a priori unknown) solution curve of \eqref{F}; (ii) starting from that point computation the solution curve $(u(s), \lambda(s))\in F^{-1}(0)$ where $s$ is the arc length; (iii) recognition and detection  of the turning point $(u^*,\lambda^*)$. To implement (ii) commonly used the following arguments.  If  $(u(s_i),\lambda(s_i))$ is a regular point (that is Rank $[D_uF|D_\lambda F](u(s_i),\lambda(s_i))=n$) then the curve $(u(s), \lambda(s))$ exists at least locally on some open interval around $s_i$ by the Implicit Function Theorem.  On this local interval along the curve the following differential equation
\begin{equation}\label{curav}
D_uF(u(s), \lambda(s)) \dot{u}(s)+D_\lambda F(u(s), \lambda(s))\dot{\lambda}(s)=0
\end{equation}
has to be satisfied. In most of continuation methods see e.g. \cite{doedK, keller1,  seydel, watson} \eqref{curav} is used to find the predictor $\tilde{u}(s_{i+1}')=u(s_i)+\alpha\cdot \dot{u}(s_i)$, where $\alpha$ is chosen in such a way that it allows then to apply Newton like iterations to find the correction $u(s_{i+1})$. Although this method is successfully applied to a various classes of problems and is an active area of researches one can not be said that it has no disadvantage see e.g. \cite{seydel}. There are certain difficulties in applying this method to large dimensional non-linear systems arising in the spatial discretization of partial differential equations.  These methods require a suitable choice of the initial point $u(s_0)$ in the neighborhood of the unknown solution curve $(u(s), \lambda(s))$, as in the indirect method, or sufficiently close to (a priori unknown) turning point in the direct methods see e.g. \cite{seydel}. In practice, the initial point actually guessed or selected from an a priori analysis.
In the case where initial point is chosen far from the unknown turning point the computation by these methods become time consuming. 
There are some difficulties with the step length control problem.
To our knowledge, there is no general strategy or the geometrical understanding in choosing of the length of the predictor step, i.e. the length of $\alpha$. 

In paper \cite{ilFunc}, it has been proposed a method where the existence of turning points of nonlinear partial differential equations is established by finding critical values of the so-called extended functional.  In particular, to find turning points for the equations of the form \eqref{F} it has been suggested the minimax  principle of the following type
\begin{equation}\label{P}
\tag{$\mathcal{P}$}
\lambda^*=\sup_{u \in S} \inf_{\psi \in \Sigma} \frac{\left\langle T(u),\psi\right\rangle }{\left\langle G(u),\psi\right\rangle }	
\end{equation}
where $S$ and $\Sigma$ some subsets in the corresponding space of solutions of \eqref{F}. This  minimax principle has been used to solve various  theoretical problems from the theory of nonlinear partial differential equations see e.g. \cite{ilcras,ilCVE,lub}  including problems which are not related with finding turning points see e.g. \cite{acker, IlD, ilCVE}. Furthermore, when $G(u)=u$ and $T(u)=Au$, where $A$ is a nonnegative matrix, \eqref{P} coincides with well-known the Collatz-Wielandt formula  for the finding the Perron-Frobenius eigenvalue see e.g. \cite{BP}.

It is our goal in the current work to develop investigations beginning in \cite{IvanIlya}, where it has been shown that the minimax formula \eqref{P} allows, in principle, to find numerically the critical value $\lambda^*$ and the corresponding turning point. 

First, we prove some new general results justifying the use of \eqref{P} for the finding  turning point. The main results in this part are Lemma \ref{L1M} and Theorems \ref{lem1} \ref{Th2} where sufficient and necessary conditions for the existence of maximal turning point of \eqref{F} in a given domain $S \subset \mathbb{R}^n$ are obtained. The second part of the work is devoted to a construction of numerical algorithm for the finding turning points by \eqref{P}. 
It should be noted that \eqref{P} belongs to a class of nonsmooth optimization problems.
Namely, we are dealing with the maximization of the following piecewise smooth function
\begin{equation}
\lambda(u):=\inf_{\psi \in \Sigma} \frac{\left\langle T(u),\psi\right\rangle }{\left\langle G(u),\psi\right\rangle }.
\end{equation} 
Nonsmooth optimization problems  have been extensively investigated in the literature over the past several decades and there are various numerical methods for such problems (see e.g., \cite{bagir1, clarke, Demyan, DemyanSPP, kiwiel,  lemar2, rock, shor, wolfe} and the references quoted in them). Our numerical approach to \eqref{P} (see also \cite{IvanIlya}) is based on the development   of the steepest ascent method for piecewise smooth function introducing by Demyanov \cite{Demyan}. 
By this method  (see also below and \cite{IvanIlya}) the steepest ascent direction $d(u)$ of $\lambda(u)$ at $u\in S$ is determined by a quadratic programming problem. 

The main feature of the presented work is that instead of steepest ascent direction we introduce the so-called quasi-direction of steepest ascent $y(u)$ of $\lambda(u)$. In this way, the direction $y(u)$ is determined by solving  a system of linear equations which has the same structure as \eqref{curav} but in general has a smaller dimension. First of all, such a replacement enables us to simplify in a certain sense the algorithm. 
However, the main goal of this approach lies in the fact that this allows us  to more precisely perform a comparative analysis of our approach with the continuation methods. 

The paper is organized as follows. In Section 2, we recall some preliminary facts about turning point type bifurcations. In Section 3, we shortly introduce the extended functional method. 
Section 4 is devoted to the  steepest ascent direction of $\lambda(u)$. In Section 5, we prove the main theoretical result on sufficient conditions providing the existence of maximal  turning point of \eqref{F}. 
In Section 6, we introduce a general algorithm for the finding maximizing point of $\lambda(u)$ in $S$ which is based on steepest ascent direction method. In Section 7, we introduce a quasi-direction of steepest ascent and a  corresponding algorithm.
Section 8 deals with numerical experiments. Finally, Section 9 is devoted to the conclusion remarks and comparative analysis.

\section{Preliminaries}

  We call $(u^*,\lambda^*)$ the {\it turning point} of \eqref{F} if there is a $C^1$-map 
 \begin{equation}\label{defb}
(-a,a) \ni s \longmapsto (u(s),\lambda(s)) \in \mathbb{R}^n\times \mathbb{R},
\end{equation}
for some $a>0$ such that {\it
 \begin{description}
	\item[(1)] $(u(s),\lambda(s))$ satisfies to \eqref{F}  for all $s \in (-a,a)$,\\ and $(u(0),\lambda(0))=(u^*,\lambda^*)$
	\item[(2)] $\frac{d}{ds}\lambda(0)=0$,
\item[(3)] one of the following  
$$
\lambda(s) \in (-\infty,\lambda^*]~~\mbox{or} ~~\lambda(s) \in  [\lambda^*,+\infty)~~ \forall s \in (-a,a), 
$$ 
is satisfied.
\end{description}} 
\noindent
We will use the notation \textit{turning point in a wide sense} for the point $(u^*,\lambda^*)$ which satisfies {\bf (1)-(2)} with  a $C^1$-map \eqref{defb}.

Denote $d u(0)/ds=\phi^*$.  \eqref{curav} yields that $D_uF(u^*, \lambda^*)\phi^*=0$, i.e.  $\phi^* \in {\rm Ker}(D_uF(u^*, \lambda^*))$.
In the literature on numerical methods see e.g. \cite{keller1, keller2, seydel}, commonly a point $(u^*,\lambda^*)$ is said to be turning point (fold bifurcation point, simple limit point)  if 
\begin{eqnarray*}
&& {\rm \bf{a)}}~ F(u^*, \lambda^*)=0,~~ \rm{ \bf {b)}}~ {\rm dim}\, Ker(D_uF(u^*, \lambda^*))=1,\\
&&{\rm\bf{c)}}~D_\lambda F (u^*, \lambda^*)  \notin {\rm Range}\,(D_uF(u^*, \lambda^*)),~~{\rm \bf{d)}}~d^2 \lambda(0)/d s^2 \neq 0.
\end{eqnarray*}
The following statement is well known see e.g. \cite{keller1, keller2, seydel}
\begin{lem}\label{kell}
Assume \textrm{\bf a)-d) } are satisfied. Then $(u^*,\lambda^*)$ is a turning point of \eqref{F}, i.e. it satisfies {\bf (1)-(3)}.
\end{lem}
Observe that \textbf{b)} holds if and only if  ${\rm dim\, Ker}(D_uF^T(u^*, \lambda^*))=1$. Therefore \textbf{b)} can be replace by:
\par\noindent \textbf{b')}  $\exists$ $\psi^* \in \mathbb{R}^n $ such that ${\rm  Ker}(D_uF^T(u^*, \lambda^*))$ = span $\{\psi^*\}$. 
\par\noindent
Furthermore,  assumption \textbf{c)} can be rewritten in the following equivalent form 
\begin{equation}
\rm{\textbf{c')}}~	\left\langle D_\lambda F (u^*, \lambda^*),\psi^*\right\rangle\neq 0~~\mbox{that is} ~~	\left\langle G(u^*),\psi^*\right\rangle\neq 0.
\end{equation}
Here and subsequently, $\left\langle x, y\right\rangle$ denotes the scalar product in $\mathbb{R}^n$, and $||x||=\left\langle x, x\right\rangle^{1/2}$ for $x,y \in \mathbb{R}^n$. 
In the future, in the case   ${\rm  Ker}(D_uF^T(u^*, \lambda^*))$ = span $\{\psi^*\}$, we shall sometimes denote the turning point of \eqref{F} by the triple $(u^*,\psi^*, \lambda^*)$.
  
\begin{rem}
In the literature, condition \textrm{\bf d)} sometime replaced by the others that prevent 
$(u^*,\lambda^*)$ from being a hysteresis point see e.g. \cite{seydel}. In our approach we also replace it (see below Lemma \ref{lem1}).
\end{rem}

\section{Extended functional method}

In this section, we shortly introduce the extended functional method  \cite{ilFunc} in its finite dimensional setting and prove some preliminary results.

By  {\it extended functional} corresponding to \eqref{F} we mean the following map
$$
Q(u, \psi, \lambda):=\left\langle  F(u, \lambda), \psi\right\rangle,~~(u,\psi, \lambda) \in \mathbb{R}^n\times \mathbb{R}^n \times \mathbb{R}.
$$
The general idea   of the extended functional method \cite{ilFunc} consists in searching of points $(u^*,\psi^*, \lambda^*) \in \mathbb{R}^n\times \mathbb{R}^n \times \mathbb{R}$ such that $(u^*,\psi^*)$ is a critical point of $Q(u, \psi, \lambda^*)$, i.e.  
\begin{equation}\label{EQ}
	\begin{cases}
	 D_\psi Q(u^*, \psi^*, \lambda^*)=0, \\
	 D_u Q(u^*, \psi^*, \lambda^*)=0, 
	\end{cases} 
\end{equation}
where $\lambda^*$ is called a critical value. 
Note that this system is nothing more than the system 
\begin{equation}\label{DEQ}
	\begin{cases}
	 F(u^*, \lambda^*)=0, \\
	 D_u F^T(u^*, \lambda^*)(\psi^*)=0. 
	\end{cases}
\end{equation}
Thus, if the critical value $\lambda^*$ is known then $(u^*,\psi^*)$ can be found as a critical point of $Q(u, \psi, \lambda^*)$ or as a solution of \eqref{EQ} (\eqref{DEQ}).  

\begin{rem}\label{DirM}
It turns out that almost the same idea is used in practice but under assumption that the value $\lambda^*$ is known approximately.  Indeed, let us replace the second equation in system  \eqref{DEQ} by $D_u F(u^*, \lambda^*)(\phi^*)=0$ and  add, for instance,  equation $\left\langle \phi^*,r \right\rangle=1$, with some $r \in \mathbb{R}^n$. Then we obtain the so-called  Seydel's  or branching system  \cite{seydel_PhD, seydel}
\begin{equation}\label{DEQ1}
	\begin{cases}
	 F(u^*, \lambda^*)=0, \\
	 D_u F(u^*, \lambda^*)(\phi^*)=0, \\
	 \left\langle \phi^*,r \right\rangle=1,
	\end{cases}
\end{equation}
which contains  $2n+1$ equations for $2n+1$ unknown variables $(u, \phi, \lambda)\in \mathbb{R}^n\times \mathbb{R}^n \times \mathbb{R}$. This system can be solved, 
for example by Newton methods provided that an initial point $(u_0,\psi_0,\lambda_0)$ of the iteration process is chosen approximately close to (a priori unknown) point	$(u^*, \psi^*, \lambda^*)$. In literature see e.g. \cite{keener, moore, seydel_PhD, seydel}, this is called a direct method of calculating bifurcation points. 
\end{rem}

In the case when the so called zero level surface $Q(u, \psi, \lambda)=0$ is solvable with respect to $\lambda$  (see \cite{ilFunc}), i.e. it is defined by a mapping $\lambda=\Lambda(u,\psi)$, the searching of points $(u^*,\psi^*, \lambda^*)$ can be carried by finding  critical points of function $Q(u, \psi,\Lambda(u,\psi))$. 

Let $\Sigma=\overline{R}^+\setminus 0$, where $R^+$ is an open orthant of the Euclidean space $\mathbb{R}^n$, i.e. 
$$
 R^+=\{\psi=\sum_{i=1}^{n} \chi_i e_i \in \mathbb{R}^n: \chi_i > 0,~i=1,...,n\},
$$
and $\overline{R}^+$  is the closure of $R^+$ in $\mathbb{R}^n$.  Let $S$ be an open subset of 
$\mathbb{R}^n$ such that 
\begin{equation}\label{GS}
\left\langle  G(u),\psi\right\rangle >0~~\mbox{for any}~~u \in S~~\mbox{and}~~\psi \in \Sigma.	
\end{equation}
In this case we are able to introduce 
\begin{equation}\label{sadd}
	\Lambda(u,\psi):= \frac{\left\langle  T(u),\psi\right\rangle}{\left\langle  G(u),\psi\right\rangle},~~u \in S, ~\psi\in \Sigma,
\end{equation}
that is under assumption \eqref{GS} the zero level surface $Q(u, \psi, \lambda)=0$ is solvable with respect to $\lambda$. 
\begin{defi}
We call $(u^*,\lambda^*)$	a maximal (minimal) turning point of \eqref{F} in $S$ if it is a turning point of \eqref{F} and 
\begin{itemize}
	\item $u^* \in S$
	\item $\lambda^*\geq \lambda_0$ ($\lambda^*\leq \lambda_0$) for any turning point $(u_0,\lambda_0)$ of \eqref{F} such that $u_0 \in S$.
\end{itemize}
\end{defi}
In \cite{ilFunc}, to find turning point of  equations of the form \eqref{F} the following variational principles have been introduced
\begin{equation}\label{MiMa}
\lambda^*=\sup_{u \in S} \inf_{\psi\in \Sigma} \frac{\left\langle  T(u),\psi\right\rangle }{\left\langle G(u),\psi\right\rangle },	
\end{equation}
\begin{equation}\label{MaMi}
\lambda_*=\inf_{u \in S} \sup_{\psi\in \Sigma} \frac{\left\langle T(u),\psi\right\rangle }{\left\langle G(u),\psi\right\rangle }.	
\end{equation}
Henceforward,  we restrict our investigation only on the maximin formula  
 \eqref{MiMa}; the minimax formula \eqref{MaMi}  can be treated in a similar fashion.

Note that $\lambda^*>-\infty$  if $S\neq \emptyset$, $\Sigma \neq \emptyset$. Furthermore, it is not hard to prove (see \cite{ilFunc})
\begin{prop}\label{NE}
Assume $\lambda^*<+\infty$. Then  equation \eqref{F} has no solutions in $S$ for any $\lambda>\lambda^*$. 	
\end{prop}
We say that $(u^*,\psi^*) \in S\times \Sigma$ is a {\it stationary point}  of $\Lambda(u,\psi)$  if 
\begin{equation}\label{EL}
\begin{cases}
D_\psi\Lambda(u^*,\psi^*)=0,\\
D_u\Lambda(u^*,\psi^*)=0.
\end{cases}
\end{equation}
Note that \eqref{EL} is equivalent to \eqref{EQ} with $\lambda^*=	\Lambda(u^*,\psi^*)$ since
\begin{align*}
&D_\psi \Lambda(u^*,\psi^*)=\frac{1}{\left<G(u^*),\psi^*\right>} F(u^*,\lambda^*), \\
&D_u\Lambda(u^*,\psi^*)=\frac{1}{\left<G(u^*),\psi^*\right>}D_u F(u^*, \lambda^*)(\psi^*).
\end{align*}
Hence, if $(u^*,\psi^*) \in S\times \Sigma$ satisfies to \eqref{EL}, then 
$$
D_{u \psi}\Lambda(u^*,\psi^*)=\frac{1}{\left<G(u^*),\psi^*\right>}D_u F(u^*, \lambda^*), 
$$ 
and, therefore  ${\rm dim\, Ker}D_u F(u^*, \lambda^*)= {\rm dim\, Ker} D_{u \psi}\Lambda(u^*,\psi^*)$.
\begin{defi}
A point  $(u^*,\psi^*) \in S\times \Sigma$ is said to be a solution of \eqref{MiMa} if it is a 	stationary point of \, $\Lambda(u,\psi)$ and\,  $\lambda^*=\Lambda(u^*,\psi^*)$. If, in addition, ${\rm dim\, Ker}D_{u \psi} \Lambda(u^*,\psi^*)=1$ then we call it a simple solution of \eqref{MiMa}.
\end{defi}

\begin{thm}\label{lem1}
Assume $\lambda^*<+\infty$  and there exists   a simple solution $(u^*,\psi^*)$ of \eqref{MiMa}. 
Then $(u^*,\lambda^*)$ is a maximal  turning point of \eqref{F} in $S$ and \, \, ${\rm Ker}(D_uF^T(u^*, \lambda^*))$= span $\{\psi^*\}$.
\end{thm}
{\it Proof.}\, 
Since $(u^*,\psi^*)$ is a simple stationary point of $\Lambda(u,\psi)$, then  
 $F(u^*,\lambda^*)=0$, $D_uF^T(u^*, \lambda^*)(\psi^*)=0$ and  ${\rm Ker}(D_uF^T(u^*, \lambda^*))$= span $\{\psi^*\}$. Thus conditions {\bf a)}, {\bf b')} are satisfied. 
By \eqref{GS} we have $\left\langle G(u^*),\psi^*\right\rangle \neq 0$. This implies (see above) that  $D_\lambda F (u^*, \lambda^*)  \notin {\rm Range}\,(D_uF(u^*, \lambda^*))$. 
Thus  {\bf c)} satisfies also. From here and since $F(u^*,\lambda^*)=0$, the Implicit Function Theorem yields that there exists map \eqref{defb} which satisfies {\bf (1)-(2)}. Thus we have  
$$
\left\langle T(u(s)),\psi\right\rangle  -\lambda(s)\left\langle G(u(s)),\psi\right\rangle =0,~~\psi \in \mathbb{R}^n, ~~s\in (-a,a).
$$
In particular, this implies that  for any $s\in (-a,a)$:
$$
\lambda(s)=\inf_{\psi\in \Sigma} \frac{\left\langle T(u(s)),\psi\right\rangle }{\left\langle G(u(s)),\psi\right\rangle }
$$
This and \eqref{MiMa} yield that $\lambda(s) \leq \lambda^*$ for all $s\in (-a,a)$. Thus {\bf(3)} satisfies and therefore $(u^*,\lambda^*)$ is a turning point of \eqref{F}. Proposition \ref{NE} yields that this is a maximal  turning point of \eqref{F} in $S$.
\hspace*{\fill}\rule{3mm}{3mm}\\

Using similar arguments  we have 
\begin{cor}\label{cor_Sad}
Let $(u^*,\psi^*)$ be a  stationary point  of $\Lambda(u,\psi)$ in  $S \times \Sigma$ such that ${\rm dim\, Ker}D_{u \psi} \Lambda(u^*,\psi^*)=1$.
Then $(u^*,\lambda^*)$ is a turning point  of \eqref{F} in a wide sense.
\end{cor}
Let $~u \in S$. Introduce
\begin{equation}\label{eq:lamb}
	\lambda(u):=\inf_{\psi\in \Sigma} \Lambda(u,\psi)\equiv \inf_{\psi\in \Sigma} \frac{\left\langle T(u),\psi\right\rangle }{\left\langle G(u),\psi\right\rangle }.
\end{equation}
Then $ -\infty \leq \lambda(u)<+\infty$ on $S$. Since $\Lambda(u,\psi)$ is a zero homogeneous with respect to  $\psi$, then $\inf_{\psi\in \Sigma} \Lambda(u,\psi)=\inf_{\psi\in \Sigma\cap \partial B_1} \Lambda(u,\psi)$, where $\partial B_1=\{x \in \mathbb{R}^n:~||x||=1\}$. Since  $\overline{\Sigma\cap \partial B_1}$ is a compact and $T,G: \mathbb{R}^n \to \mathbb{R}^n$ are continuous functions, then $\lambda(u)$ is a continuous function on $S$. 

From now on we make the assumptions: 
\begin{description}
	\item[(H)] For any $u_0\in S$   
\begin{equation*}
S(u_0):=\{u \in S: \lambda(u)	\geq \lambda(u_0)\}~~\mbox{is bounded and}~~ \overline{S}(u_0)\subset S.
\end{equation*} 
\end{description}
Here $\overline{S}(u_0)$ is the closure of $S(u_0)$.
\begin{lem}\label{lem2}
Suppose that hypothesis   \rm{({\bf H})} is satisfied. Then 
\begin{description}
	\item[1)] $-\infty<\lambda^*<+\infty$
	\item[2)] there exists  $u^* \in S$ such that $\lambda(u^*)=\lambda^*$.
	\end{description}
\end{lem}
{\it Proof.}\,
The problem \eqref{MiMa} can be rewritten in the following form
\begin{equation}\label{mima}
	\lambda^*=\sup\{\lambda(u): ~u \in S\}.
\end{equation}
Let $(u_m) \subset S$ be a maximizer sequence of \eqref{mima}, i.e. $\lambda(u_m) \to \lambda^*$ as $m \to \infty$.
Then \rm{({\bf H})} yields that $u_m \in \overline{S}(u_1)$ for $m=1,2,...$. Note that $\overline{S}(u_1)$ is compact in $\mathbb{R}^n$ since it is bounded. Hence, there is a subsequence $u_{m_i}$, $i=1,2,...$ such that $ u_{m_i} \to u^*$ as $i\to \infty$. Then $u^* \in S$ since by \rm{({\bf H})} we have  $\overline{S}(u_1)\subset S$. The continuity of  $\lambda(u)$ implies
that $\lambda(u_{m_i}) \to \lambda(u^*)$ as $i\to \infty$. Hence,  since $\lambda(u_m) \to \lambda^*$ as $m \to \infty$ we get that $\lambda^*=\lambda(u^*)$ and $\lambda^*<\infty$. 
\hspace*{\fill}\rule{3mm}{3mm}\\

\section{The  steepest ascent direction}

Observe that minimization problem  \eqref{eq:lamb}
can be rewritten as the following linear programming problem  
\begin{eqnarray}
&& {\rm minimize}~~~\left\langle  T(u),\psi \right\rangle  , \label{lin}\\ 
&&{\rm subject} ~{\rm to}~~~\left\langle  G(u),\psi \right\rangle  =1,~~
\psi \in \Sigma. \label{simpl}
\end{eqnarray}
From this  we see that $\lambda(u)$ is bounded and the minimum in \eqref{lin}-\eqref{simpl} is achieved at one of the vertices  of  polygon \eqref{simpl}. Thus \eqref{eq:lamb} is equivalent to
\begin{eqnarray}
\lambda(u)=\min_{i} \frac{\left\langle  T(u),{e}_i\right\rangle }{\left\langle  G(u), {e}_i\right\rangle }\equiv \min_{i} f_i(u), ~~u \in S,
\label{L}
\end{eqnarray}
where  $ f_i(u)=\left\langle  T(u),{e}_i\right\rangle /\left\langle  G(u), {e}_i\right\rangle $, $i=1,2,...,n$.  Thus, minimax problem \eqref{MiMa} can be replaced by the following
\begin{equation}\label{MiMaS}
	\lambda^*=\max_{u\in S} \min_{i} f_i(u).
\end{equation}
By the assumption, $T,G: \mathbb{R}^n \to \mathbb{R}^n$ are continuously differentiable functions. Therefore $f_i \in C^1(S)$, $i=1,2,...,n$ since by \eqref{GS} one has $\left\langle  G(u),e_i\right\rangle > 0$ in $S$. This implies that $\lambda(u)$ is a piecewise continuously differentiable function on $S$. Furthermore, (see e.g. \cite{Demyan} and \cite{IvanIlya})
$\lambda(u)$ is directionally differentiable  in $S$  with respect to any vector $d \in \mathbb{R}^n$ and the directional derivative is defined by
\begin{equation}\label{pl1}
\lambda'(u;d)=\min_{i}\{\left\langle \nabla f_i(u), d\right\rangle: i \in N(u)\} ,
\end{equation}
where $ \nabla f_i(u)=(\frac{\partial f_i(u)}{\partial u_1},...,\frac{\partial f_i(u)}{\partial u_n})^T$ and  
$$
N(u)=\{i\in [1:n]:~f_i(u)=\lambda(u)\}.
$$
Hereinafter, $N:=|N(u)|$ denotes the number of elements in $N(u)$, and the series $i_1,...,i_{N} \in N(u)$ denotes the arrangement of the set $N(u)$  in ascending sequence $i_1<i_2<\ldots<i_{N}$.
Note that $u$ satisfies \eqref{F} if and only if the equalities $f_i(u)=\lambda(u)$  hold for all $i=1,2,...,n$ or the same when   $|N(u)|=n$.

Following \cite{Demyan} we call a maximizer $\hat{d}(u) \in \mathbb{R}^n$ of 
\begin{equation}\label{pl2}
\hat{\sigma}(u):=\lambda'(u;\hat{d}(u))=\max\{\lambda'(u;d):  ||d||=1\}.	
\end{equation}
a \textit{direction of steepest ascent} of $\lambda(u)$ at $u\in S$ if $\hat{\sigma}(u) >0$. In this case we denote
$$
\nabla \lambda(u):=\lambda'(u;\hat{d}(u))\cdot\hat{d}(u)\equiv \hat{\sigma}(u)\cdot\hat{d}(u)
$$
Consider the subdifferential of $\lambda(u)$ at $u \in S$ 
\begin{align}\label{CH1}
	\partial \lambda(u):&={\rm conv}\{\nabla f_i(u): i\in N(u)\}
\end{align}
where for any set $C$, conv$C$ denotes the convex hull of $C$.
In the sequel, Nr $C$ denotes the point of smallest Euclidean norm in the conv$C$, i.e. Nr$C$
is a nearest point  from the origin $0_n$ to conv$C$. By Demyanov-Malozemov's Theorem (see Theorem 3.3 in \cite{Demyan}, see also \cite{DemyanSPP}, \cite{IvanIlya}) one has
\begin{equation}\label{DM}
\nabla \lambda(u)={\rm Nr}\, \partial \lambda(u),	
\end{equation}
and 
\begin{equation}\label{DM1}
\hat{d}(u)=\nabla \lambda(u)/||\nabla \lambda(u)||,~~~ \hat{\sigma}(u)=	||\nabla \lambda(u)||.
\end{equation}
Observe that
\begin{equation}\label{dif}
\nabla f_i(u)=\frac{1}{ G_i(u)}\left[\nabla T_i(u) -\lambda(u)\nabla G_i(u)\right], ~~i=1,2,..., n.
\end{equation}
Let ${\cal M}$ be an arbitrary subset of $ \{[1:n]\}$, $M:=|{\cal M}|$ and $i_1,...,i_{M} \in {\cal M}$ is an arrangement of the set ${\cal M}$  in ascending sequence $i_1<i_2<\ldots<i_M$. Introduce matrices
\begin{equation}\label{A}
A(u)=\left(G_i(u)\cdot \nabla f_i(u)\right)_{ 1\leq i\leq n}^T, ~~A_{\cal M}=\left(G_{i_k}(u)\cdot \nabla f_{i_k}(u)\right)_{ 1\leq k\leq M}^T,	
\end{equation}
\begin{equation}\label{A2}
{\cal A}(u)=\left(\nabla f_i(u)\right)_{1\leq i\leq n}^T,~~{\cal A}_{\cal M}=\left(\nabla f_{i_k}(u)\right)_{ 1\leq k\leq M}^T,	
\end{equation}
and
$$
\Gamma_{{\cal M}}:={\cal A}_{\cal M}^T{\cal A}_{\cal M}.
$$
Then by \eqref{DM}, \eqref{DM1} maximization problem \eqref{pl2} is equivalent to the following 
quadratic programming problem 
\begin{align}\label{SMin}
\hat{\sigma}^2(u)=\min_{\alpha}\{&\alpha^T \Gamma_{N(u)} \alpha:~\alpha \in \mathbb{R}^N,~\\
&\sum_{i=1}^{N} \alpha_i=1,~~\alpha_i\geq 0,~ i=1,...,N\},\nonumber
\end{align}
so that if $\hat{\alpha}(u)$ is a minimizer of \eqref{SMin}, then
\begin{equation}\label{NL}
\nabla \lambda(u)=\sum_{k=1}^N \hat{\alpha}_k \nabla f_{i_k}(u)\equiv {\cal A}_{N(u)} \hat{\alpha}(u),	
\end{equation}
and 
\begin{equation}\label{NL1}
\hat{d}(u)=\nabla \lambda(u)/\hat{\sigma}(u)\equiv \frac{{\cal A}_{N(u)} \hat{\alpha}(u)}{||{\cal A}_{N(u)} \hat{\alpha}(u)||}
\end{equation} 
is a maximizer of \eqref{pl2} (see also \cite{IvanIlya}). Note that the minimizer $\hat{\alpha}(u)$  of \eqref{SMin} always exists. 

Denote 
$$
 \Sigma_{N(u)}=\overline{\{\psi=\Sigma_{i=1}^{N} \alpha_{k_i} e_i \in \mathbb{R}^N: \alpha_i > 0,~i=1,...,N \}}\setminus 0.
$$

\begin{cor}\label{cor:SN}
$\hat{\sigma}(u)=0$ if and only if there exists $\psi \in \Sigma_{N(u)}$ such that ${\cal A}_{N(u)}\psi=0$, i.e.  $\psi \in {\rm Ker}(D_uF^T(u, \lambda))$ with $\lambda=\Lambda(u,\psi)$.  
\end{cor}
{\it Proof.}\,
Assume $\hat{\sigma}(u)=0$. Then there exists a minimizer $\psi \in \Sigma_{N(u)}$ of \eqref{SMin}  such that
$$
0=\left\langle \Gamma_{N(u)} \psi ,\psi \right\rangle=\left\langle {\cal A}_{N(u)}\psi ,{\cal A}_{N(u)}\psi \right\rangle\equiv ||{\cal A}_{N(u)}\psi ||^2.
$$ 
But this is possible only if ${\cal A}_{N(u)}\psi=0$. The proof of the inverse statement is trivial.
\hspace*{\fill}\rule{3mm}{3mm}\\

By \cite{Demyan} (see Theorem 2.1 in \cite{Demyan}) for any $u\in S$ and $d \in \mathbb{R}^n$ , $||d||=1$ one has
\begin{equation}\label{Dh}
\lambda(u+\tau d)=\lambda(u)+\tau \lambda'(u;d)+\bar{o}(d;\tau)	
\end{equation}
for sufficiently small $\tau \in \mathbb{R}$, where $\bar{o}(d;\tau)/\tau \to 0$ as $\tau \to 0$.

From this and by \eqref{pl1}, \eqref{NL}, \eqref{NL1} we have 
\begin{cor}\label{cor:sid2}
Assume that $\hat{\sigma}(u)>0$. Then there exist $\tau_0>0$ such that
$$
\lambda(u+\tau \hat{d}(u))>\lambda(u)
$$
for any $\tau \in (0,\tau_0)$, where $\hat{d}(u)$ is given by \eqref{NL1}
\end{cor}

\begin{lem}\label{L1}
 Let  $u^* \in S$ be  a maximizer of $\lambda(u)$, i.e. 	$\lambda(u^*)=\max_{u\in S}\lambda(u)$. Then 
$\hat{\sigma}(u^*)=0$ and there exists $\psi^* \in \Sigma_{N(u^*)}$ such that ${\cal A}_{N(u^*)}(u^*)\psi^*=0$, i.e.  $\psi \in {\rm Ker}(D_uF^T(u^*, \lambda^*))$ with $\lambda^*=\Lambda(u^*,\psi^*)$.
 \end{lem}
{\it Proof.}\,
 Suppose, contrary to our claim, that $\hat{\sigma}(u^*)>0$. Then by Corollary  \ref{cor:sid2} 
there is a steepest ascent direction $\hat{d}(u^*)$  such that $\lambda(u^*+\tau \hat{d}(u^*))>\lambda(u^*)$ for sufficiently small $\tau>0$.
Evidently, $u^*+\tau \hat{d}(u^*) \in S$ for small  $\tau>0$ since  $S$ is an open set. 
 Thus we get a contradiction and consequently $\hat{\sigma}(u^*)=0$. The proof of the last part of the lemma follows from Corollary \ref{cor:SN}.
\hspace*{\fill}\rule{3mm}{3mm}\\ 
 

\begin{rem}
The condition $\hat{\sigma}(u)>0$ does not imply that the matrix $  A_{N(u)}(u)$ is nonsingular. Indeed, it is possible that there is $\psi_0 \not\in \Sigma_{N(u)}$ such that $A_{N(u)}\psi_0=0$  while $A_{N(u)}\psi\neq 0$ for any $\psi \in\Sigma_{N(u)}$.	
\end{rem}

From the above, we have also the following necessary condition for a maximal (minimal) turning point of \eqref{F} in a given $S \subset \mathbb{R}^n$

\begin{lem}\label{L1M} 
Let $S \subset \mathbb{R}^n$ be a given open subset of 
$\mathbb{R}^n$ such that \eqref{GS} holds. Assume $(u^*,\lambda^*)$	is a maximal (minimal) turning point of \eqref{F} in $S$. Then  $\hat{\sigma}(u^*)=0$ and there exists $\psi^* \in {\rm Ker}(D_uF^T(u^*, \lambda^*))$ such that  $\psi^*_i\geq 0$, $i=1,...,n$.
 \end{lem}
{\it Proof.}\,
Since \eqref{GS} holds we are able to introduce the function 
$$
\lambda(u)=\inf_{\psi\in \Sigma} \frac{\left\langle T(u),\psi\right\rangle }{\left\langle G(u),\psi\right\rangle },
$$
and consider the minimization problem \eqref{SMin}. Suppose, contrary to our claim, that $\hat{\sigma}(u^*)>0$. Then arguing as above in the proof of Lemma \ref{L1} we see that this is impossible. Thus  $\hat{\sigma}(u^*)=0$ and the proof follows from Corollary \ref{cor:SN}.
\hspace*{\fill}\rule{3mm}{3mm}\\

\section{Existence of the maximal turning point}

In this section, we derive some sufficient conditions when minimax problem \eqref{MiMaS} gives a turning point of \eqref{F}.

From now on we make, in addition, the following assumption:
\begin{description}
	\item[(R)]   Rank ${\cal A}(u) \geq n-1$ for any  $u \in S$.
\end{description}

\begin{thm}\label{Th2}
Assume $T,G: \mathbb{R}^n \to \mathbb{R}^n$ are continuously differentiable functions. 
Suppose that hypothesis  \eqref{GS},  {\rm ({\bf H})}, {\rm ({\bf R})} are satisfied. Then there exists a solution $(u^*,\psi^*)$ of \eqref{MiMa} such that $(u^*,\lambda^*)$ is a maximal turning point  of \eqref{F} in $S$. Furthermore, $\psi^* \in {\rm Ker}(D_uF^T(u^*, \lambda^*))$ such that  $\psi^*_i> 0$, $i=1,...,n$.
\end{thm}
{\it Proof.}\,
By Lemma \ref{lem2} there exists $u^* \in S$ such that $\lambda^*=\lambda(u^*)$. Consequently Lemma \ref{L1} yields 
$\hat{\sigma}(u^*)=0$, and therefore  there is $\psi^* \in \Sigma_{N(u^*)}$ such that  ${\cal A}_{N(u^*)} \psi^* =0$. By assumption \rm{({\bf R})} this is possible only if $|N(u^*)|=n$. This yields that $\lambda^*=f_i(u^*)$ for all $i=1,...,n$,
that is $u^*$ satisfies \eqref{F}. Furthermore, the vector $\psi^* \in {\rm Ker}({\cal A}(u^*))$ is defined uniquely  up to scalar multiplication since  Rank ${\cal A}(u^*) \geq n-1$. Thus,  $(u^*,\psi^*)$ is a simple solution of \eqref{MiMaS}. This implies by Theorem \ref{lem1} that  $(u^*,\lambda^*)$ is a maximal turning point  of \eqref{F} in $S$. 
\hspace*{\fill}\rule{3mm}{3mm}\\

\begin{rem}
Assumption {\rm ({\bf R})} is not so restrictive. For instance, it is commonly satisfied for systems arising in the spatial discretization of partial differential equations (see below). 
\end{rem}

\begin{rem} 
From the proof of Theorem \ref{Th2} it can be seen that the existence of solution $(u^*,\psi^*)$ of \eqref{MiMaS} will be still hold if we replace {\rm ({\bf R})} by the following 
\begin{description}
	\item[(RW)]   dim ${\rm Ker}({\cal A}(u))\cap \Sigma\leq 1$ for   $u \in S$.
\end{description}
\end{rem}
Under assumption {\rm ({\bf R})} we can strengthen Corollary \ref{cor_Sad} as follows
\begin{cor}\label{cor:Sad2}
Assume  {\rm ({\bf R})} holds. 
Let $(u^*,\psi^*)$ be a stationary point  of $\Lambda(u,\psi)$ in  $S \times \Sigma$.
Then $(u^*,\lambda^*)$ is a turning point of \eqref{F} in a wide sense.
\end{cor}

\section{A general algorithm}

In this section, we discuss how to construct an algorithm for the finding turning point of \eqref{F}  based on the above theory.

Henceforth we always assume that hypothesis {\rm ({\bf H})} and {\rm ({\bf R})} are satisfied. 
We say $u^*_\varepsilon$ is a solution of \eqref{F} with accuracy $\varepsilon>0$, if
\begin{equation}
|f_i(u^*_\varepsilon)-\lambda(u^*_\varepsilon)|<\varepsilon~~\mbox{for all}~~i=1,2,...,n. 	
\end{equation}
Let us denote $N_\varepsilon(u)=\{i\in [1:n]:~f_i(u)-\lambda(u)< \varepsilon\}$ for $u \in S$. Thus $u$ is a solution of \eqref{F} with  accuracy $\varepsilon>0$ if and only if $|N_\varepsilon(u)|=n$.

From above we know  that $\hat{\sigma}(u^*)=0$ is the necessary condition for $(u^*,\lambda^*)$ to be a maximal turning point of \eqref{F}. We will detect this condition up to accuracy $\sigma>0$. 
\begin{defi}
Let $\varepsilon>0$, $\sigma>0$. We  call $(u^*_{(\varepsilon,\sigma)},\psi^*_{(\varepsilon,\sigma)},\lambda^*_{(\varepsilon,\sigma)})$ the {$(\varepsilon,\sigma)$-maximal turning point of \eqref{F}}  (or shortly call the $(\varepsilon,\sigma)$-m.turning point) if 
\begin{description}
	\item[(i)] $\lambda^*_{(\varepsilon,\sigma)}=\lambda(u^*_{(\varepsilon,\sigma)})$, $|N_\varepsilon(u^*_{(\varepsilon,\sigma)})|=n$, 
	\item[(ii)] $\psi^*_{(\varepsilon,\sigma)}$ is a minimizer of
		\begin{align}\label{SMinES}
\sigma^2_{N_\varepsilon(u^*_{(\varepsilon,\sigma)})}=\min_{\psi}\{&\psi^T \Gamma_{N(u^*_{(\varepsilon,\sigma)})} \psi:~\psi \in \mathbb{R}^N,~\\
&\sum_{i=1}^{N} \psi_i=1,~~\psi_i\geq 0,~ i=1,...,N\},\nonumber
\end{align}
		\item[(ii)] 		$\sigma_{N_\varepsilon(u^*_{(\varepsilon,\sigma)})}<\sigma$.
\end{description}
\end{defi}

With respect to the above it can be proposed the following general algorithm for the finding 
turning point of \eqref{F}.  
\par\noindent
{\bf ALGORITHM 1  (Algorithm of the steepest ascent direction)}. 
\par\noindent
{\it 
Choose an initial point $u_0 \in S$ and accuracies $\varepsilon>0$, $\sigma>0$. 

 For $k:=0, 1,2,...$ until $(\varepsilon,\sigma)$-m.turning point $(u^*_{(\varepsilon,\sigma)},\psi^*_{(\varepsilon,\sigma)})$ is found
\begin{description}
	\item[1)] Find 
	$$
	\lambda(u^k)=\min_{1\leq i\leq n} f_i(u^k),
	$$ 
	input 
	$$
	N_\varepsilon(u^k)=\{i\in [1:n]:~|f_i(u^k)-\lambda(u^k)|<\varepsilon\},~~N=|N_\varepsilon(u^k)|
	$$
	\item[2)] Find  $\psi^k$ by minimization 
\begin{align}
\sigma_{N_\varepsilon(u^k)}^2=\min_{\psi}\{&\psi^T \Gamma_{N_\varepsilon(u^k)} \psi:~\sum_{i=1}^{N} \psi_i=1,~~\psi_i\geq 0,~ i=1,...,N\}.\label{Kmin}
\end{align}

\item[3)] If $\sigma_{N_\varepsilon(u^k)}<\sigma$ and $|N_\varepsilon(u^k)|=n$, then go to Step 4), else 	
\begin{description}
	\item[3.1)]  introduce the steepest ascent direction
$$
d^k=	\frac{{\cal A}(u^k) \psi^k}{|| {\cal A}(u^k) \psi^k||}
$$
	\item[3.2)] 
	 find step length $\tau^k$ by  golden section search rule applying to 
	$$
	\kappa(\tau):=\lambda(u^k+\tau d^k):=\min_{1\leq i\leq n} \{f_i(u^k+\tau d^k)\},~\tau\geq 0,~~u^k+\tau d^k \in S.
	$$
	\item[3.3)] introduce $u^{k+1}:=u^k+\tau^k d^k$ and return to Step 1).
	\end{description}
	
	\item[4)] Output the $(\varepsilon,\sigma)$-m.turning point: $u^*_{(\varepsilon,\sigma)}:=u^k$, $\psi^*_{(\varepsilon,\sigma)}:=\psi^k$, $\lambda^*_{(\varepsilon,\sigma)}:=\lambda(u^k)$.
\end{description}
}

This algorithm has been implemented in \cite{IvanIlya}. A justification of the applicability of the algorithm for the finding turning point follows from the next lemmas

\begin{lem}\label{lem:alg1}
Assume that {\rm ({\bf H})}, {\rm ({\bf R})} are satisfied. Then for any given 	$\varepsilon>0$, $\sigma>0$ 
there exists $k= k(\varepsilon,\sigma)>0$ such that $|N_\varepsilon(u^k)|=n$ and $\sigma_{N_\varepsilon(u^k)}<\sigma$.
\end{lem} 

\begin{lem}\label{lem:exis} 
Assume that conditions {\rm ({\bf H})}, {\rm ({\bf R})} are satisfied.\\	Let $(u^*_{(\varepsilon,\sigma)},\psi^*_{(\varepsilon,\sigma)}, \lambda^*_{(\varepsilon,\sigma)})_{\varepsilon>0,\sigma>0}$ be a set of $(\varepsilon,\sigma)$-m.turning points of \eqref{F}. Then there exists a limit point $(u^*,\psi^*,\lambda^*)$  of  $(u^*_{(\varepsilon,\sigma)},\psi^*_{(\varepsilon,\sigma)},\lambda(u^*_{(\varepsilon,\sigma)}))$ \\  as $\varepsilon, \sigma \to 0$ such that $(u^*,\lambda^*)$ is a turning point of \eqref{F} in a wide sense and ${\rm  Ker}(D_uF^T(u^*, \lambda^*))$ = span $\{\psi^*\}$.  
\end{lem}

The proofs of these lemmas is similar to given below proofs of Lemmas \ref{lem:alg2}, \ref{lem:exis2}, respectively.


\section{A quasi-direction of steepest ascent}

 Analysis of Algorithm 1 shows that the most costly step is {\bf 2)}, i.e. the finding the steepest ascent direction. 
To implement this step it is necessary to find minimizer of \eqref{Kmin}. Note that \eqref{Kmin} is a quadratic programming problem. The theory of the numerical solution of such problems is well developed see e.g. \cite{flet, vasil} and to solve \eqref{Kmin} one of the method from this theory can be used.
In our paper \cite{IvanIlya}, the numerical implementation of Algorithm 1 was based precisely on this idea, i.e. on the iterative finding the minimizer of quadratic programming problem \eqref{Kmin}.

In this section, instead of the steepest ascent direction we introduce its approximation, the so-called quasi-direction of steepest ascent. This direction can be found by solving a system of linear equations that is less costly than minimization of quadratic programming problem. 
However, the main goal of such a replacement consists in the fact that the system of linear equations of the quasi-direction of steepest ascent method are similar to that used in the continuation methods \cite{doedK}, \cite{keller1}, \cite{seydel}. Thus, we will be able to compare more precisely the  two approaches,  the extended functional and continuation methods.

Let $u \in S$ and   $\hat{\alpha}:=\hat{\alpha}(u) \in \mathbb{R}^N$ be a corresponding minimizer of \eqref{SMin}. Then by  Karush-Kuhn-Tucker conditions there exist constants $\mu_0$, $\mu_{i}$, $i=1,..., N\equiv |N(u)|$ such that 
\begin{equation}
\label{KKT}
\left\{
\begin{array}{l}
\Gamma_{N(u)}\hat{\alpha}-\mu_0 1_{N}-\sum_{i=1}^N \mu_{k_i} {e}_i=0 \\ \\
\left\langle  1_{N},\hat{\alpha} \right\rangle  \equiv \sum_{i=1}^N \hat{\alpha}_i=1\\ \\
\hat{\alpha}_i \geq 0,~\mu_i\geq 0,~\mu_i\hat{\alpha}_i =0,~~i=1,2,...,N.  
\end{array}
\right.
\end{equation}
Here $\Gamma_{N(u)}:={\cal A}_{N(u)}^T{\cal A}_{N(u)}$, $1_{N}=(1,...,1)^T \in \mathbb{R}^{N}$.
Thus, the  minimizer $\hat{\alpha}(u)$ can be obtained by solving \eqref{KKT}. Note that the dimension of this problem is
$2N(u)+1$ with unknown variables $\alpha_i$, $\mu_i  \in \mathbb{R}$, $i=1,... N$ and $\mu_0 \in \mathbb{R}$ . 

In our approach instead of this we will use the following less dimensional  system 
\begin{equation}
\label{D}
\begin{cases}
 \Gamma_{N(u)} \alpha=\delta \cdot 1_{N}, \\[1em]
  \sum\limits_{i=1}^{N}\alpha_i = 1.
\end{cases}
\end{equation}
where $\alpha=(\alpha_1,...,\alpha_{N})^T$ and $\delta \in \mathbb{R}$. 
\begin{rem}\label{RemC}
In the case $N(u)=n$, $\delta=1$ if one put   $v={\cal A}_{N(u)} \alpha$, the first equation in \eqref{D} is easily converted see \cite{IvanIlya} to the Davidenko-Abbott system see  \cite{abbot,david}
\begin{equation}\label{ceq3}
	D_uF(u,\lambda) v= -D_\lambda F(u,\lambda),
\end{equation}
which lies at the core of the continuation methods see e.g. \cite{doedK, keller1,  seydel}. 
\end{rem}
Observe that by \eqref{D} we have
\begin{equation}\label{gamm}
\left\langle \Gamma_{N(u)} \alpha,\alpha\right\rangle =|| \sum_{k=1}^N \nabla f_{i_k}(u)\alpha_k||^2=\delta  \geq 0.	
\end{equation}
Furthermore, hypothesis  {\rm ({\bf R})} implies that for 
any $u \in S$  system \eqref{D} has a unique solution $(\alpha(u), \delta)$. 

Assume that $|N(u)|>1$. Let $i_1,...,i_{N} \in N(u)$ be an arrangement of the set $N(u)$  such that $i_1<i_2<\ldots<i_{N}$.  Consider the following affine space
$$
{\cal L}_{N(u)}=
\{v=\sum_{i\in N(u)} \beta_i \nabla f_{i}(u):~\sum_{i\in N(u)}\beta_i=1\}. 
$$
Introduce
$$
Y(u)=\sum_{k=1}^N \nabla f_{i_k}(u)\alpha_k(u),
$$
where $\alpha(u)$ satisfies \eqref{D}. 
Assume that $\delta>0$. Then by  (\ref{D}) we have 
\begin{equation}\label{deleq}
	\langle  Y(u),  \nabla f_{i}(u)\rangle=\langle  {\cal A}_{N(u)}\alpha,  \nabla f_{i}(u)\rangle =\langle  \Gamma_{N(u)}\alpha, {e}_i \rangle\equiv \delta,
\end{equation}
$\forall i \in N(u)$. Hence 
$\langle  Y(u),  \nabla f_{i}(u)-\nabla f_{j}(u)\rangle=0$, $\forall i,j \in N(u)$. From 
{\rm ({\bf R})} and since $\delta>0$ it is easy to infer that $\nabla f_{i}(u)\neq \nabla f_{i}(u)$ for all $i,j \in N(u)$, $i\neq j$. 
This shows that $Y(u)$ is an orthogonal vector to ${\cal L}_{N(u)}$. In other words, 
$Y(u)$ is a nearest point  from the origin $0_n$ to the affine space  ${\cal L}_{N(u)}$. 
Note that the  subdifferential $\partial \lambda(u)$ lies on ${\cal L}_{N(u)}$.

Introduce $y(u)=Y(u)/||Y(u)||$. We call vector $y(u)$ the {\it quasi-direction of steepest ascent}.  Recall that $\nabla \lambda(u)$ is a nearest point  from the origin $0_n$ to the subdifferential $\partial \lambda(u)$.

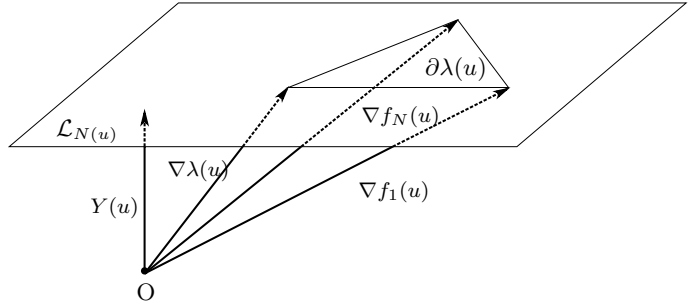
\begin{figure}
\centering

  \setlength{\unitlength}{1bp}%
  \begin{picture}(257.67, 116.55)(0,0)
    \put(0,0){\includegraphics{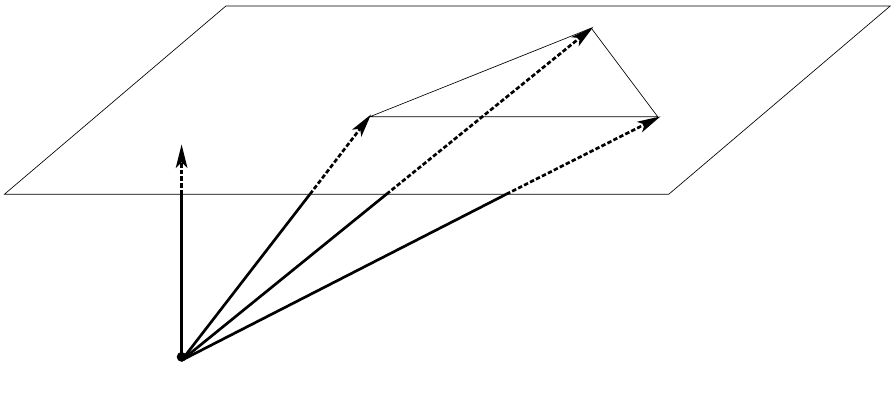}}
  \put(19.03,65.20){\fontsize{9}{10.76}\selectfont ${\cal L}_{N(u)}$}
  \put(31.86,35.89){\fontsize{8}{9.22}\selectfont $Y(u)$}
  \put(60.88,50.10){\fontsize{8}{9.22}\selectfont $\nabla \lambda(u)$}
  \put(158.06,88.12){\fontsize{9}{10.76}\selectfont $\partial \lambda(u)$}
  \put(132.97,42.02){\fontsize{8}{9.22}\selectfont $\nabla f_1(u)$}
  \put(134.08,70.51){\fontsize{8}{9.22}\selectfont $\nabla f_N(u)$}
  \put(49.19,3.20){\fontsize{9}{10.76}\selectfont O}
  \end{picture} 
  \caption{Quasi-direction of steepest ascent $Y(u)$ and direction of steepest ascent $\nabla \lambda(u)$} 
 \label{fig:f1} 
\end{figure}

\begin{lem}\label{invSim}
Suppose that hypothesis   {\rm ({\bf R})} is satisfied. Let $u\in S$ and $(\alpha,\delta)$ be a solution of (\ref{D}) such that $\delta> 0$. Then
\begin{description}
	\item[a)]  $Y(u)=\nabla \lambda(u)$ if and only if $\alpha_k>0$, $\forall k=1,..., |N(u)|$. Furthermore, if $\alpha_k>0$, $\forall k=1,..., |N(u)|$, then
	$\delta=\sigma^2_{N(u)}$.
	\item[b)] If there are subsets $N_1(u), N_2(u)$ such that $N_1(u)\cup N_2(u) = N(u)$ and   $\alpha_k>0$, $\forall k \in N_1(u)$,  whereas $\alpha_k\leq 0$, $\forall k \in N_2(u)$, then  $\nabla \lambda(u)$ lies on the boundary  
		\begin{equation*}\label{CH}
	\partial_{N_1(u)}\lambda(u):=\{\sum_{i\in {N_1(u)}} \nabla f_i(u)\zeta_i:\sum_{i\in{N_1(u)}}\zeta_i=1,~~\zeta_i\geq 0, ~i\in {N_1(u)}\}.
\end{equation*}
of $\partial \lambda(u)$,  i.e. $\nabla \lambda(u) \in \partial_{N_1(u)}\lambda(u)$. 
\end{description}
\end{lem}
{\it Proof.}\,
The proof of {\bf a)} is evident. Assume $\alpha$  satisfies  {\bf b)}.
Then  $Y(u)$ belongs to the set
\begin{eqnarray*}
C:=\{w=\sum_{i\in N(u)} \beta_i\nabla f_{i}(u)&: &~\sum_{i\in N(u)}\beta_i =1,~\\  
\beta_i >0,&& \forall i \in N_1(u), ~~ \beta_j \leq 0,~~ \forall j \in N_2(u) \}.
\end{eqnarray*}
The intersection of $C$ with $\partial \lambda(u)$ coincides with $\partial_{N_1(u)}\lambda(u)$. Hence, the nearest point $Z'$  from $Y(u)$ to $\partial \lambda(u)$
belongs to $\partial_{N_1(u)}\lambda(u)$. But evidently $Z'=\nabla \lambda(u)$ (see Fig. \ref{fig:f1}).  
\hspace*{\fill}\rule{3mm}{3mm}\\

\begin{lem}\label{zeroSim}
Suppose that hypothesis   {\rm ({\bf R})} is satisfied. Let $u\in S$ and $(\alpha,\delta)$ be a solution of (\ref{D}). 
\begin{description}
	\item[a)] If $\delta=0$, then $|N(u)|=n$ and $\alpha_k\neq 0$, $\forall k=1,..., n$.
	\item[b)] If $\delta=0$ and  $\alpha_k>0$, $\forall k=1,...,n$, then $\hat{\sigma}(u)=0$, and there exists $\psi \in \Sigma$ such that  ${\cal A}\psi=0$, i.e.  $\psi \in {\rm Ker}(D_uF^T(u, \lambda))$ with $\lambda=\Lambda(u,\psi)$. 
	\item[c)] If $\delta=0$ and there are subsets $N_1(u), N_2(u)$ such that $N_1(u)\cup N_2(u) = N(u)$ and  $\alpha_k>0$, $\forall k \in N_1(u)$,  whereas $\alpha_k\leq 0$, $\forall k \in N_2(u)$, then $\nabla \lambda(u)$ lies on the boundary  $\partial_{N_1(u)}\lambda(u)$ of $\partial \lambda(u)$,  i.e. $\nabla \lambda(u) \in \partial_{N_1(u)}\lambda(u) $.
\end{description}
\end{lem}
{\it Proof.}\,
Statement {\bf a)} is a direct consequence of the assumption ({\bf R}). The proof of {\bf b)} follows from \eqref{SMin}, \eqref{gamm} and Corollary \ref{cor:SN}. To prove {\bf c)} consider the set 
\begin{eqnarray*}
Q:=\{w=\sum_{i\in N(u)} \beta_i w_i&|&~\sum_{i\in N(u)}\beta_i =1,~~ \beta_i\geq 0,~~ 
w_i=\nabla f_{i}(u),\\&&~~ \forall i \in N_1(u),~ w_i=-\nabla f_{i}(u),~~ \forall i \in N_2(u) \}.
\end{eqnarray*}
Then $0_n \in Q$ and $Q\cap \partial_{N(u)}= \partial_{N_1(u)}$. This implies that the nearest point $\nabla \lambda(u)$ from $0_n$ to  $\partial\lambda(u)$ lies on  $\partial_{N_1(u)}\lambda(u)$. 
\hspace*{\fill}\rule{3mm}{3mm}\\

Lemmas \ref{invSim}, \ref{zeroSim} allow us to build  an algorithm for the finding of the  steepest ascent direction of $\lambda(u)$ by solving (\ref{D}). However, in present paper we are not going to use the steepest ascent direction $d(u)$, since  we  use the  quasi-direction of steepest ascent $y(u)$ of $\lambda(u)$.

Let $u \in S$. Introduce the following bordered matrix cf. \cite{keller1}
\begin{equation}
	{\cal M}_{N(u)}=\left( \begin{array}{cc}
	\Gamma_{N(u)}& -1_{N(u)} \\
	1_{N(u)}^T & 0  \\
	 \end{array} \right).
\end{equation}
Denote
\begin{eqnarray*}
&&t:=t(u)=\left(\begin{array}{c}
\alpha(u) \\
\delta(u) \\
 \end{array}\right) \in \mathbb{R}^{N+1},\\
&&\alpha:=\alpha(u) \in \mathbb{R}^N, \delta:=\delta(u) \in \mathbb{R},\\~~~
 &&q_N=\left(\begin{array}{c}
0_N \\
1 \\
 \end{array}\right) \in \mathbb{R}^{N+1},
 \end{eqnarray*}
Then (\ref{D}) is equivalent to 
\begin{equation}\label{bordE}
	{\cal M}_{N(u)} t=q_{N(u)}.
\end{equation}

\begin{defi}
Let $\varepsilon>0$, $\delta>0$. We  call $(u^*_{(\varepsilon,\delta)},\psi^*_{(\varepsilon,\delta)},\lambda^*_{(\varepsilon,\delta)})$ the $(\varepsilon,\delta)$-turning point of \eqref{F}  if 
\begin{description}
	\item[(i)] $\lambda^*_{(\varepsilon,\delta)}=\lambda(u^*_{(\varepsilon,\delta)})$, $|N_\varepsilon(u^*_{(\varepsilon,\delta)})|=n$, 
	\item[(ii)] $(\psi^*_{(\varepsilon,\delta)}, \delta_{N_\varepsilon(u^*_{(\varepsilon,\delta)})})$ satisfies 
		\begin{align}\label{SMinESq}
{\cal M}_{N(u^*_{(\varepsilon,\delta)})} \left(\begin{array}{c}
\psi^*_{(\varepsilon,\delta)}\\
\delta_{N_\varepsilon(u^*_{(\varepsilon,\delta)}} \\
 \end{array}\right)=q_{N(u^*_{(\varepsilon,\delta)})}
\end{align}
		\item[(ii)] $\delta_{N_\varepsilon(u^*_{(\varepsilon,\delta)})}<\delta$.		
\end{description}
\end{defi}

We propose the following algorithm of the quasi-direction of steepest ascent for the finding  the $(\varepsilon,\delta)$-turning point $(u^*_{(\varepsilon,\delta)},\psi^*_{(\varepsilon,\delta)},\lambda^*_{(\varepsilon,\delta)})$ of \eqref{F}  

{\bf ALGORITHM 2. (AQDSA)}. 
\par\noindent
{\it 
Choose an initial point $u_0 \in S$ and accuracies $\varepsilon>0$, $\delta>0$. 

 For $k:=0, 1,2,...$ until $(\varepsilon,\delta)$-turning point $(u^*_{(\varepsilon,\delta)},\psi^*_{(\varepsilon,\delta)},\lambda^*_{(\varepsilon,\delta)})$ is found:
\begin{description}
	\item[1)] Find $$
	\lambda(u^k)=\min_{1\leq i\leq n} f_i(u^k),
	$$
	and input 
	$$
	N_\varepsilon(u^k)=\{i\in [1:n]:~|f_i(u^k)-\lambda(u^k)|<\varepsilon\},~~N=|N_\varepsilon(u^k)|.
	$$
	\item[2)] Find  $\delta^k$ and $\alpha^k$ by solving 
	$$
	{\cal M}_{N_\varepsilon(u^k)} t^k=q_N~~\mbox{where}~~t^k=\left(\begin{array}{c}
\alpha^k \\
\delta^k \\
 \end{array}\right) ,~~q_N=\left(\begin{array}{c}
0 \\
1 \\
 \end{array}\right) \in \mathbb{R}^{N+1}.
	$$

\item[3)] If $\delta^k<\delta$ and $N=n$, then go to Step 4), otherwise 	
\begin{description}
	\item[3.1)]  Introduce the quasi-direction of steepest ascent
$$
y^k=Y^k/||Y^k||~~\mbox{where}~~Y^k:=\sum_{j=1}^N \nabla f_{i_j}(u^k) \alpha_j.
$$
	\item[3.2)] find step length $\tau^k$ by golden section search rule  applying to 
	$$
	\kappa(\tau):=\lambda(u^k+\tau y^k):=\min_{1\leq i\leq n} \{f_i(u^k+\tau y^k)\},~\tau\geq 0,~~u^k+\tau d^k \in S.
	$$
	\item[3.3)] introduce $u^{k+1}:=u^k+\tau^k y^k$ and return to Step 1).
	\end{description}
	
	\item[4)] Output the $(\varepsilon,\delta)$-turning point: $u^*_{(\varepsilon,\delta)}:=u^k$, $\psi^*_{(\varepsilon,\delta)}:=\alpha^k$, $\lambda^*_{(\varepsilon,\delta)}:=\lambda(u^k)$.
\end{description}
}

Let us show that this algorithm gives indeed the $(\varepsilon, \delta)$-turning point.

First we prove
\begin{prop}\label{prop:R}
Let $u \in S$ and $\varepsilon >0$. Then there is $r:=r(u, \varepsilon)>0$ such that
\begin{equation}
N_\varepsilon(u)=N_\varepsilon(v) ~~~\mbox{for any}~~~ v \in B_r(u):=\{v:~||u-v||<r\}.	
\end{equation}
\end{prop}
{\it Proof.}\,
By definition $f_i(u)-\lambda(u)< \varepsilon$ for any $i \in N_\varepsilon(u)$. Since the functions $f_i(u)$ and $\lambda(u)$ are continuous then there is a neighborhood $B_{r_0}(u)$ of $u$ with some $r_0>0$ such that the inequalities 
$f_i(v)-\lambda(v)< \varepsilon$ will hold for any $v \in B_{r_0}(u)$. Hence $N_\varepsilon(u) \subseteq N_\varepsilon(v)$ for $v \in B_{r_0}(u)$. Clear  that  $N_\varepsilon(u) = N_{\varepsilon_0}(u)$ for some $\varepsilon_0>\varepsilon$. Let $B_{r_1}(u)$ such that   
$|f_i(u)-f_i(v)|< (1/2)(\varepsilon_0-\varepsilon)$, $i=1,2,...,n$ and  $|\lambda(u)-\lambda(v)|< (1/2)(\varepsilon_0-\varepsilon)$ for  
any $v \in  B_{r_1}(u)$. Then 
\begin{eqnarray*}
|f_i(u)-\lambda(u)| \leqslant |f_i(u)-f_i(v)| + |f_i(v)-\lambda(v)|+|\lambda(v)-\lambda(u)|<\varepsilon_0
\end{eqnarray*} 
for any $v \in  B_{r_1}(u)$ and $i \in N_\varepsilon(v)$. Thus $N_\varepsilon(v)\subseteq N_{\varepsilon_0}(u)=N_\varepsilon(u)$. Hence $N_\varepsilon(u) = N_\varepsilon(v)$ for $v \in B_r(u)$, where $r=\min\{r_0,r_1\}$.
\hspace*{\fill}\rule{3mm}{3mm}\\

\begin{prop}\label{prop:Inc2}
 Let $(u^k)$ be a iteration sequence defined by Algorithm 2. Suppose that $\delta^k>0$. Then 
\begin{equation}\label{vozr2}
\lambda(u_{k+1})>\lambda(u^k).	
\end{equation}
\end{prop}
{\it Proof.}\,
By Proposition \ref{prop:R} one has $N_\varepsilon(u^k+\tau y^k)=N_\varepsilon(u^k)$ for sufficiently small $\tau>0$.
Furthermore, evidently $N(u^k+\tau y^k) \subseteq N_\varepsilon(u^k+\tau y^k)$.
Thus, for sufficiently small $\tau$ we have
\begin{align}
\lambda(u_{k+1})&=\max_{\tau\geq 0}\lambda(u^k+\tau y^k) \geq \lambda(u^k+\tau y^k)=\nonumber\\
&\min_{i \in N(u^k+\tau y^k)}[f_i(u^k)+\tau \left\langle \nabla f_i(u^k), y^k \right\rangle  
+\phi_i(\tau,u^k)] \geq \nonumber\\
&\min_{i \in N_\varepsilon(u^k+\tau y^k)}[f_i(u^k)+\tau \left\langle \nabla f_i(u^k), y^k \right\rangle  
+\phi_i(\tau,u^k)] \geq \nonumber  \\ 
&\lambda(u^k)+\tau \min_{i \in N_\varepsilon(u^k)}\left\langle\nabla f_i(u^k), y^k \right\rangle  +\min_{i \in N_\varepsilon(u_{k})}\phi_i(\tau,u^k) \label{eq:voz22}
\end{align}
where $\phi_i(\tau,u^k)=o(\tau)$, $i=1,2,...,n$ as $\tau \to 0$. By (\ref{D}) (see also \eqref{deleq}) we have 
$$
\langle  \nabla f_{i}(u^k), y^k  \rangle=\langle  \nabla f_{i}(u^k), Y^k  \rangle/||Y(u^k)||=\delta^k /||Y(u^k)||.
$$
This implies that
\begin{equation}\label{se2}
	\min_{i \in N_\varepsilon(u^k)}\left\langle\nabla f_i(u^k), y^k \right\rangle =\delta^k /||Y(u^k)||>0.
\end{equation}
Therefore the sum of the last two terms in \eqref{eq:voz22} is positive for sufficiently small $\tau$. This yields \eqref{vozr2}. 
\hspace*{\fill}\rule{3mm}{3mm}\\

\begin{lem}\label{lem:alg2}
Assume that \eqref{GS}, {\rm ({\bf H})}, {\rm ({\bf R})} are satisfied. Let $\varepsilon>0$, $\delta>0$ are given. Then for any starting point $u_0 \in S$ there exists a finite number 	 
$k= k(\varepsilon,\delta)>0$ such that $|N_\varepsilon(u^{k})|=n$ and $\delta^{k}<\delta$ for the iteration sequence $(u^k, \delta^k)$ defined by Algorithm 2.
\end{lem} 
{\it Proof.}\,
Let $(u^k)$ be the iteration sequence defined by Algorithm 2. Then Lemma \ref{lem2} and Proposition \ref{prop:Inc2} imply that there is a unique limit value $\hat{\lambda} = \lim^{k\to \infty}\lambda(u^k)$ such that
 $\lambda(u^k)\leq \hat{\lambda} \leq \lambda(u^*)$, $k=1,2,...,$.

Suppose, contrary to our claim, that for any $k=0,1,...,$ it holds one of the following: $|N_\varepsilon(u^k)|<n$ or/and $\delta^k\equiv \delta(u^k)>\delta$. 
By \rm{({\bf H})} $u^k \in \overline{S}(u_0)$ for $k=1,2,...$. Then the compactness of $\overline{S}(u_0)$ implies
the existence of a subsequence of $u^k$ (we denote it again
by $u^k$ ) such that  $ u^k \to \hat{u}$ as $k\to \infty$. Then $\hat{u} \in S$ since   $\overline{S}(u_0)\subset S$. Using Proposition \ref{prop:R} it is not hard to see from \eqref{bordE} that $\delta(u)$ is a continuous function on $S$. Therefore $\delta(u^k) \to \delta(\hat{u})$ as $k\to \infty$. By Proposition \ref{prop:R} it follows  that $|N_\varepsilon(u^k)| = |N_\varepsilon(\hat{u})|$ for sufficiently large $k$. From these and our assumption there are only two possibilities  
\begin{description}
	\item[a)] $|N_\varepsilon(\hat{u})|\leq n$,  $\delta(\hat{u})>\delta$,
	\item[b)] $|N_\varepsilon(\hat{u})|<n$, $\delta(\hat{u})\leq \delta$.
\end{description}
Assume that {\bf a)} holds. By Proposition \ref{prop:R} there is $r:=r(u^*, \varepsilon)>0$ such that
$|N_\varepsilon(u)| = |N_\varepsilon(\hat{u})|$ for any $u \in B_r:=B_r(\hat{u})$. Obviously, there is $\tau_0>0$ such that
$u+\tau h \in B_r$ for any $u\in B_{r/2}(\hat{u})$, $h \in \mathbb{R}^n: ||h||=1$ and $\tau \in (0,\tau_0)$. Then $|N_\varepsilon(u+\tau h)| = |N_\varepsilon(\hat{u})|$ for these $u,h,\tau$. Since $ u^k \to \hat{u}$  and $\delta(u^k) \to \delta(\hat{u})$ as $k\to \infty$, then  we can find a number $K>0$ such that $u^k \in B_{r/2}$ and $\delta(u^k)>\delta$ for every $k>K$.

However, as in \eqref{eq:voz22} for $k>K$ we have 
\begin{eqnarray}
\lambda(u^{k+1})&\geq& \lambda(u^k)+\tau \min_{i \in N_\varepsilon(u^k)}\left\langle \nabla f_i(u^k), y^k \right\rangle  + \label{eq:voz2}\\
&&\min_{i \in N(u^{k})}\phi_i(\tau,u^k) >\hat{\lambda} + \omega^k +\tau \delta + \min_{i \in N(u^{k})}\phi_i(\tau,u^k) \nonumber 
\end{eqnarray}
where $ \omega^k=(\lambda(u^k) -\hat{\lambda})$. The continuously differentiability of functions $f_i$  yields that $\sup_{u \in  B_{r/2}}|\phi_i(\tau,u)|/\tau \to 0$ as $\tau \to 0$. Consequently there exists $\tau_1 \in (0,\tau_0)$ such that
$$
\tau \delta + \min_{i \in N(u^{k})}\phi_i(\tau,u^k)>\frac{1}{2}\tau \delta
$$
for any $\tau \in (0,\tau_1)$ and $k>K$. Now taking into account that $\omega^k \to 0$ as $k\to \infty$ we obtain from \eqref{eq:voz2} that $\lambda(u^{k+1})>\hat{\lambda}$ for sufficiently large $k$. However this contradicts to the inequalities $\lambda(u^k)\leq \hat{\lambda}$, $k=1,2,...$. Thus {\bf a)} can not to be satisfied.  

Suppose that {\bf b)} holds. Then $\delta(\hat{u})=0$. Indeed, if it is not hold then we can repeat the above arguments that has been used in the case {\bf a)} and obtain again the contradiction. However, $\delta(\hat{u})=0$ implies by \eqref{bordE} that ${\cal A}(u)$ is a singular matrix. But this is impossible under assumptions {\rm ({\bf R})} and $|N_\varepsilon(\hat{u})|<n$. This completes the proof.
\hspace*{\fill}\rule{3mm}{3mm}\\

From this we have
\begin{cor}\label{lem:alg2c}
Assume that \eqref{GS}, {\rm ({\bf H})}, {\rm ({\bf R})} are satisfied. Let $\varepsilon>0$, $\delta>0$ are given. Then for any starting point $u_0 \in S$ Algorithm 2 gives the $(\varepsilon,\delta)$-turning point $(u^*_{(\varepsilon,\delta)},\psi^*_{(\varepsilon,\delta)},\lambda^*_{(\varepsilon,\delta)})$ in finite steps.
\end{cor}
Let us now prove
\begin{lem}\label{lem:exis2} 
Assume that conditions {\rm ({\bf H})}, {\rm ({\bf R})} are satisfied.\\	Let $(u^*_{(\varepsilon,\delta)},\alpha^*_{(\varepsilon,\delta)}, \lambda(u^*_{(\varepsilon,\delta)}))_{\varepsilon>0,\delta>0}$ be a set of $(\varepsilon,\delta)$-turning points of \eqref{F}. Then there exists a limit point $(u^*,\alpha^*,\lambda^*)$  of  $(u^*_{(\varepsilon,\delta)},\alpha^*_{(\varepsilon,\delta)},\lambda(u^*_{(\varepsilon,\delta)}))$  as  $\varepsilon, \delta \to 0$ such that $(u^*,\lambda^*)$ is a turning point of \eqref{F} in a wide sense and ${\rm  Ker}(D_uF^T(u^*, \lambda^*))$ = span $\{\alpha^*\}$.
\end{lem}
{\it Proof.}\,
From assumption {\rm ({\bf H})} using the same arguments as in the proof of Lemma \ref{lem2} it can be shown that there is a subsequence $u^*_i:=u^*_{(\varepsilon_i,\delta_i)}$ with $\varepsilon_i, \delta_i \to 0$ such that  $u^*_i \to u^*$  as $i \to \infty$. Then the continuity of $\lambda(u)$ implies that there is a limit value $\lambda^*=\lim_{i\to \infty}\lambda(u^*_i)$. Since $|N_{\varepsilon}(u^*_{(\varepsilon,\delta)})|=n$, then evidently $|N(u^*)|=n$. This yields that $u^*$ satisfies \eqref{F}. 
Furthermore, we have $\delta(u^*_i) \to \delta(u^*)$ as $i \to \infty$ and  $\delta(u^*_i)<\delta_i$ for $i=1,2,...$. This implies that $\delta(u^*)=0$ since $\delta_i \to 0$.  By  \eqref{bordE} this is possible only if ${\cal A}(u)$ is a singular matrix. Furthermore, evidently there is a limit $\alpha^*_{(\varepsilon_i,\delta_i)} \to \alpha^*$ as $i \to \infty$.  Therefore passing to the limit in the equality $\Gamma (u^*_i)\alpha^*_{(\varepsilon_i,\delta_i)}=\delta(u^*_i)1_n$  one get $\Gamma (u^*)\alpha^*=0$. Thus  $\alpha^* \in {\rm  Ker}(D_uF^T(u^*, \lambda^*))$.  Now, taking into account assumption {\rm ({\bf R})} we obtain the proof.
\hspace*{\fill}\rule{3mm}{3mm}\\

One should keep in mind that the point  $(u^*,\alpha^*)$ obtained by Lemma \ref{lem:exis2} as a limit of the set of 
$(\varepsilon,\delta)$-turning points $(u^*_{(\varepsilon,\delta)},\alpha^*_{(\varepsilon,\delta)})$ does not necessary satisfy to the condition 
$\sigma(u^*)=0$, since it is possible $\alpha^*_i< 0$ for some $i \in [1:n]$ (see Lemma \ref{zeroSim}, {\bf c)}). However, if $\sigma(u^*)>0$ then by Lemma \ref{L1}  $(u^*,\alpha^*)$ can not be maximal turning point of \eqref{F} in $S$ and will  not be detected by Algorithm 1. 
On the other hand, $(u^*,\alpha^*)$ may indeed be a turning point, since $\delta=0$ and consequently ${\rm  Ker}(D_uF^T(u^*, \lambda^*))\neq \emptyset$.   Summarizing, one can say that, in general, Algorithm 1 is more efficient in the finding of the maximal turning point of \eqref{F} in $S$, while Algorithm 2 may be useful in searching for  all type of turning points of \eqref{F}. 

Below in numerical experiments, we are dealing with problems where the corresponding matrices ${\cal A}(u)$ have a sparse structure, namely they are tridiagonal. In this case, the  direct application of Algorithm $2$  is time consuming. The situation can be improved if we take the value of $\varepsilon$ in Step $1$ initially big enough and then iteratively reduced it to the required accuracy. This allows us to increase the set $N_\varepsilon(u^k)$ in the initial steps, and thereby increase the number of involved variables $(u^k_j)$ varying function $\lambda(u)$. 

Thus we use the following modified algorithm of the quasi-direction of steepest ascent 

{\bf ALGORITHM 3. (MAQDSA)}. 
\par\noindent
{\it 
Choose an initial point $u_0 \in S$ and accuracies $\varepsilon>0$, $\delta>0$. 

 For $k:=0, 1,2,...$ until $(\varepsilon,\delta)$-turning point $(u^*_{(\varepsilon,\delta)},\psi^*_{(\varepsilon,\delta)},\lambda^*_{(\varepsilon,\delta)})$ is found:
\begin{description}
   \item[1)] Find 
		\begin{align*}
		\lambda(u^k):&=\min_{1\leq i\leq n} f_i(u^k),\\
		\mu(u^k):&= \max_{1\leq i\leq n} f_i(u^k),
	\end{align*}
	input $\epsilon^k:=(\mu(u^k)-\lambda(u^k))/2$
	\item[2)] If $\epsilon^k<\varepsilon$, then $\epsilon^k:=\varepsilon$. 
	\item[3)] Input 
	$$
	N(u^k)=\{i\in [1:n]:~|f_i(u^k)-\lambda(u^k)|<\epsilon^k\},~~N=|N(u^k)|.
	$$
	\item[4)] Find  $\delta^k$ and $\alpha^k$ by solving 
	$$
	{\cal M}_{N_\varepsilon(u^k)} t^k=q_N~~\mbox{where}~~t^k=\left(\begin{array}{c}
\alpha^k \\
\delta^k \\
 \end{array}\right) ,~~q_N=\left(\begin{array}{c}
0 \\
1 \\
 \end{array}\right) \in \mathbb{R}^{N+1}.
	$$

\item[5)]  If $\delta^k<\delta$ and $N=n$, then go to Step 6).

 \begin{description}
	\item[5.1)]  Introduce the quasi-direction of steepest ascent
$$
y^k=Y^k/||Y^k||~~\mbox{where}~~Y^k:=\sum_{j=1}^N \nabla f_{i_j}(u^k) \alpha_j.
$$
	\item[5.2)] find step length $\tau^k$ by golden section search rule  applying to
	$$
	\kappa(\tau):=\lambda(u^k+\tau y^k):=\min_{1\leq i\leq n} \{f_i(u^k+\tau y^k)\},~\tau\geq 0,~~u^k+\tau d^k \in S.
	$$
	\item[5.3)] introduce $u^{k+1}:=u^k+\tau^k y^k$, $\lambda(u^{k+1}):= \kappa(\tau^k)$, $\epsilon^{k+1}:=\epsilon^k$ and return to Step 2).
	\end{description}
	
	\item[6)] If $\epsilon^k>\varepsilon$, then 
	$$
	\epsilon^k:=\epsilon^k/2~~ \mbox{and go to Step}~~ 2),
	$$
	else output the $(\varepsilon,\delta)$-turning point: $u^*_{(\varepsilon,\delta)}:=u^k$, $\psi^*_{(\varepsilon,\delta)}:=\alpha^k$, $\lambda^*_{(\varepsilon,\delta)}:=\lambda(u^k)$.
\end{description}
}

\section{Numerical implementations}

In this section the results of numerical experiments performed with the modified algorithm quasi-direction of steepest ascent (MAQDSA) are presented.  The algorithm was implemented in MATLAB and run on a PC AMD Athlon II X2 of 2.71 GHz CPU and 1.75 GB of RAM. The algorithm was tested on several examples and the results for two cases: Bratu-Gelfand problem and elliptic equation with convex-concave nonlinearity are given below. In both these examples we  consider branches of positive solutions and as a set of $S$ in \eqref{MiMa}  we take an open positive orthant of the Euclidean space $\mathbb{R}^n$: 
$$
 S=\{\sum_{i=1}^{n} u_i e_i :~  u_i > 0,~i=1,...,n\}.
$$
The number of iterations (ItN) and computing time (CPU) are reported 
as a measure of the  performance. To test the performance of MAQDSA, we compare it with the performance of the  numerical continuation packages {\small \textsc{cl\_matcont}}  see e.g. \cite{matcont}. 
In order to provide a fair comparison,  we unplugged part of the functions  in the {\small \textsc{cl\_matcont}} so that it sought only limit points. To apply the {\small \textsc{cl\_matcont}}  the specification of the  initial point $(u_{\lambda_0},\lambda_0)$ on the solution path is required. An approximate of $u_{\lambda_0}$  at an arbitrary fixed value $\lambda_0$ was obtained by the standard routine from  MATLAB. Subsequently, index $_{MC}$ stands for {\small \textsc{cl\_matcont}}, e.g. $(u_{MC}, \lambda_{MC})$ denotes turning point found by {\small \textsc{cl\_matcont}}. 

We report computations  performed for $n = 100$. For the stopping criterion we used $\varepsilon = 10^{-6}$ and test different values $\delta = 10^{-9}, \delta = 10^{-10}, \delta = 10^{-11},\delta = 10^{-12}$. Furthermore, we tested the MAQDSA for different  initial points: $u_0 =  0.1 \cdot 1_n$, $u_0 = 1_n$, $u_0 = 10 \cdot 1_n$, where $1_n:=(1,...,1)^T \in \mathbb{R}^n$.

{\bf Example 1. (convex-concave problem)} 
 Consider the boundary value  problem with convex-concave nonlinearity
\begin{equation}
\label{convex_concave} 
	\begin{array}{l}
	-\Delta u = \lambda u^{q}+u^{\gamma},~ x \in \Omega, \\
	u |_{\partial \Omega}=0,
	\end{array}
\end{equation}
where $\Omega = (0,1)$  and $\partial \Omega$ denotes the boundary of $\Omega$  and $0<q<1<\gamma$. We discretized $\Omega$ by a uniform grid with grid points $x_i=i \cdot h$, $1\leqslant i\leqslant n$, where $h=1/(n+1)$. For the second derivatives at $n$ mesh points we used  a standard second-order finite difference approximation. This yields the system of $n$ nonlinear algebraic equations

\begin{equation}
\label{dis1} 
	\begin{array}{l}
	-\frac{u_{i+1}-2 u_i +u_{i-1}}{h^2}=  \lambda u_i^q+u_i^{\gamma},~ 1\leqslant i\leqslant n, \\
	u_0=u_{n+1}=0.
	\end{array}
\end{equation}
Then the functions $f_i: \mathbb{R}^n \to \mathbb{R}$ are given by 
\begin{eqnarray*}
&& f_i(u)=\frac{-u_{i+1}+2 u_i -u_{i-1}-h^2 u_i^\gamma}{h^2 u_i^q},~~i=2,...,n-1,	\\
&& f_1(u)=\frac{-u_{2}+2 u_1 -h^2 u_1^\gamma}{h^2 u_1^q},\\
&& f_n(u)=\frac{2 u_n -u_{n-1}-h^2 u_n^\gamma}{h^2 u_n^q}.
\end{eqnarray*}
By direct calculation of the corresponding matrix ${\cal A}(u)$ it can be seen that it is a tridiagonal. This implies that the matrix ${\cal A}(u)$ satisfies to condition  {\rm ({\bf R})}. In \cite{IvanIlya} it has  been justified that condition {\rm ({\bf H})} is also satisfied. Thus we may apply Theorem \ref{Th2} and therefore there exists a maximal turning point of discretized convex-concave problem  \eqref{dis1}. Furthermore, the turning point $(u^*,\psi^*,\lambda^*)$ can be found as a solution of the corresponding maximin problem \eqref{MiMaS}  applying the steepest ascent direction or quasi-direction of steepest ascent (MAQDSA) algorithm. As we know by Lemmas \ref{lem:alg2}, \ref{lem:exis2} for any given $\varepsilon>0$, $\delta>0$ and  any starting point $u_0 \in S$ MAQDSA gives in finit steps the $(\varepsilon,\delta)$-turning points of \eqref{dis1}.  

The value $\lambda^*$ of the turning point $(u^*,\lambda^*)$ depends from the parameters $q$ and $\gamma$, where $0<q<1<\gamma$. As examples, we present  two different (in a certain sense) cases:  $q=0.5, \gamma = 2$ ($\lambda_{MC}= 11.643872$), and $q=0.1, \gamma = 1.5$ ($\lambda_{MC}=93.140742$ is bigger). 

In order to select the appropriate stop criterion (comparable by accuracy with the  {\small \textsc{cl\_matcont}}), 
MAQDSA was tested at different $\delta = 10^{-9}, \delta = 10^{-10}, \delta = 10^{-11},\delta = 10^{-12}$ (see Table \ref{tab:delta}).
Distances are measured in the norms  $\left\|x\right\| = \sqrt{\sum\limits_{i=1}^{n}x_i^2} $, 
$\left\|x\right\|_\infty= \max_{1 \leqslant i \leqslant n}{\left|x_i\right|} $, $(u^*,\lambda^*)$ denotes $(\varepsilon,\delta)$-turning point obtained by MAQDSA. From Table \ref{tab:delta} we see that the appropriate stop criterion for $\delta$ can be taken $\delta = 10^{-11}$ or $\delta = 10^{-12}$. Similar results appear for other cases including below the Bratu-Gelfand problem.  

The test results for the performances of  {\small \textsc{cl\_matcont}} (first colum) and  MAQDSA (the last three columns)  with $\delta=10^{-9}$ and different $q, \gamma$  are reported in Tables \ref{tab:CC_100_5} - \ref{tab:CC_100_2}. 
 \medskip

 {\bf Example 2. (Bratu-Gelfand problem)}  Consider the Bratu-Gelfand problem
\begin{equation}
\label{Bratu} 
	\begin{array}{l}
	-\Delta u = \lambda e^u,~ x \in \Omega, \\
	u |_{\partial \Omega}=0,
	\end{array}
\end{equation}
where $\Omega$  and $\partial \Omega$ are the same as in \eqref{convex_concave}. As above in Example 1, we consider the following discretization of \eqref{Bratu}
\begin{equation}
\label{DBratu} 
	\begin{array}{l}
	-\frac{u_{i+1}-2 u_i +u_{i-1}}{h^2}= \lambda e^{u_i},~ 1\leqslant i\leqslant n, \\ 
 	u_0=u_{n+1}=0.
	\end{array}
\end{equation}
In this case the functions $f_i: \mathbb{R}^n \to \mathbb{R}$ are given by 
\begin{eqnarray*}
&& f_i(u)=\frac{-u_{i+1}+2 u_i -u_{i-1}}{h^2  e^{u_i}},~~i=2,...,n-1,	\\
&& f_1(u)=\frac{-u_{2}+2 u_1 }{h^2  e^{u_1}},~~~ f_n(u)=\frac{2 u_n -u_{n-1}}{h^2  e^{u_n}}.
\end{eqnarray*}  
Arguing as above see also \cite{IvanIlya}, it can be conclude that  MAQDSA  gives a $(\varepsilon,\delta)$-turning points of discretized Bratu-Gelfand problem \eqref{DBratu}.

The test results for the performances of {\small \textsc{cl\_matcont}}  and  MAQDSA with $\delta=10^{-9}$   are reported in Table \ref{tab:B_100}.

\begin{table}[tbp]
\caption{Problem with convex-concave nonlinearity \newline  in case $q=0.5, \gamma = 2, u_0 = 1_n,  \lambda_{MC} = 11.643872$}
\begin{tabular}{c c c c c}
\hline
 & $\delta = 10^{-9}$ & $\delta = 10^{-10}$ & $\delta = 10^{-11}$ & $\delta = 10^{-12}$ \\
\hline
	ItN & 107 & 109 & 115 & 118 \\
	$\left|\lambda_{MC} - \lambda^*\right|$ & $3.1 \cdot 10^{-8}$ & $3.1 \cdot 10^{-8}$ & $3.1 \cdot 10^{-8}$ & $3.1 \cdot 10^{-8}$ \\  
 $\left\|u_{MC} - u^*\right\|$ & $1.75 \cdot 10^{-3}$ & $8.30 \cdot 10^{-4}$ & $1.04 \cdot 10^{-4}$ & $6.09 \cdot 10^{-5}$ \\
 $\left\|u_{MC} - u^*\right\|_\infty$ & $2.49 \cdot 10^{-4}$ & $1.18 \cdot 10^{-4}$ & $1.50 \cdot 10^{-5}$ & $9.05 \cdot 10^{-6}$ \\
\hline
\end{tabular}
\label{tab:delta}
\end{table}

\begin{table}[tbp]
\caption{Problem with convex-concave nonlinearity \newline in case  $q=0.5, \gamma = 2 \, (\lambda_0=1,~~\lambda^* = 11.643872)$}
\begin{tabular}{c c c c c c}
\hline
 & {\small \textsc{cl\_matcont}} & & \multicolumn{3}{c}{MAQDSA}\\
\hhline{~-~---}
 & $ u_0 =  u_{\lambda_0}$ & & $u_0 =  0.1 \cdot 1_n$ & $ u_0 = 1_n$ & $ u_0 = 10 \cdot 1_n$ \\
\hline
 ItN & 320 &  & 111  & 107 & 106 \\  
 CPU & 1.2 &  & 0.73 & 0.75 & 0.75 \\
\hline
\end{tabular}
\label{tab:CC_100_5}
\end{table}

\begin{table}[tbp]
\caption{Problem with convex-concave nonlinearity \newline in case  $q=0.1, \gamma = 1.5 \, (\lambda_0=1,~~\lambda^* = 93.140742)$}
\begin{tabular}{c c c c c c}
\hline
 & {\small \textsc{cl\_matcont}} & & \multicolumn{3}{c}{MAQDSA} \\
\hhline{~-~---}
 & $u_0 =  u_{\lambda_0}$ & & $u_0 =  0.1 \cdot 1_n$ & $ u_0 = 1_n$ & $ u_0 = 10 \cdot 1_n$ \\
\hline
 ItN & 3600 &  & 269 & 191 & 142 \\  
 CPU & 14.7 &  & 2.16 & 1.64 & 1.31 \\
\hline
\end{tabular}
\label{tab:CC_100_2}
\end{table}
\begin{table}[tbp]
\caption{Bratu-Gelfand problem $(\lambda_0=1,~\lambda^* = 3.513652)$}
\begin{tabular}{c c c c c c}
\hline
 & {\small \textsc{cl\_matcont}} & & \multicolumn{3}{c}{MAQDSA} \\
\hhline{~-~---}
 & $ u_0 =  u_{\lambda_0}$ & & $u_0 =  0.1 \cdot 1_n$ & $ u_0 = 1_n$ & $ u_0 = 10 \cdot 1_n$ \\
\hline
 ItN & 88 &  & 121 & 101 & 136 \\  
 CPU & 0.5 &  & 0.70 & 0.70 & 0.83 \\
\hline
\end{tabular}
\label{tab:B_100}
\end{table}
The examples illustrate that the performance of MAQDSA generally comparable  with {\small \textsc{cl\_matcont}} (and sometimes even surpasses) and can be used for the finding turning points of the problems type \eqref{F}.  We remark that our algorithm has not (yet) been tuned for efficiency, and we expect that the computational effort required by carefully designed algorithm, to be smaller.

\section{Conclusion Remarks}

\begin{figure}
\centering
  \setlength{\unitlength}{1bp}%
  \begin{minipage}[h]{0.45\linewidth} 
  	\begin{picture}(154.31, 152.09)(0,0)
    	\put(0,0){\includegraphics{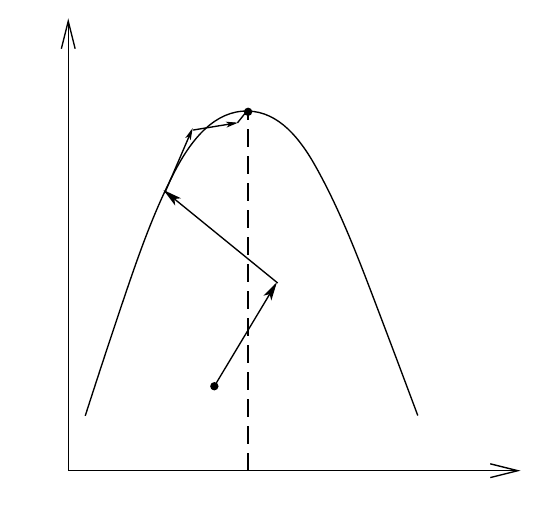}}
  	\put(137.84,7.42){\fontsize{10}{12.48}\selectfont $u$}
  	\put(70.62,7.66){\fontsize{10}{12.48}\selectfont $u^*$}
  	\put(5.18,131.50){\fontsize{10}{12.48}\selectfont $\lambda$}
  	\put(49.51,32.47){\fontsize{7}{8.54}\selectfont $\lambda(u_0)$}
  	\put(84.17,71.26){\fontsize{7}{8.54}\selectfont $\lambda(u_1)$}
  	\put(51.71,98.77){\fontsize{7}{8.54}\selectfont $\lambda(u_2)$}
  	\put(71.68,125.48){\fontsize{7}{8.54}\selectfont $\lambda(u^*)$}
  	\put(114.81,62.59){\fontsize{7}{8.54}\selectfont $\lambda(u_{\lambda}) = f_1(u_{\lambda}) = \ldots = f_n(u_{\lambda}) $}
  	\end{picture}
 		\caption{Schematic representation of the iterative procedure by the extended functional method.} 
 		\label{fig:f2} 
 	\end{minipage}
	\hfill
	\begin{minipage}[h]{0.45\linewidth}
		\begin{picture}(172.77, 141.11)(0,0)
  	\put(0,0){\includegraphics{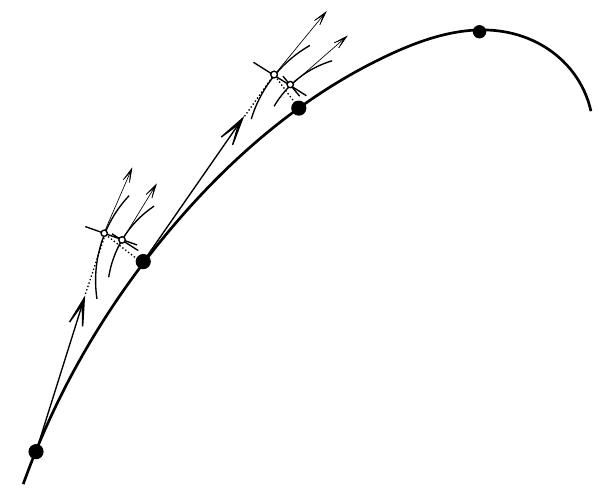}}
  	\put(15.49,11.49){\fontsize{7}{7.77}\selectfont $u(s_i)$}
  	\put(45.96,65.58){\fontsize{7}{7.77}\selectfont $u(s_{i+1})$}
  	\put(90.92,108.12){\fontsize{7}{7.77}\selectfont $u(s_{i+2})$}
  	\put(2.58,48.68){\fontsize{7}{7.77}\selectfont $\dot{u}(s_i)$}
  	\put(6.95,79.96){\fontsize{7}{7.77}\selectfont $\tilde{u}(s_{i+1}')$}
  	\put(40.83,104.00){\fontsize{7}{7.77}\selectfont $\dot{u}(s_{i+1})$}
  	\put(133.49,124.60){\fontsize{7}{7.77}\selectfont $u^*$}
  	\end{picture}
  	\caption{Schematic representation of the iterative procedure by the continuation method.} 
 		\label{fig:matcont}
 	\end{minipage}	
 \end{figure}

In this paper, we  continued the elaboration of new paradigm of the finding bifurcations turning point type launched in \cite{IvanIlya}. The main distinguishing feature of this approach is that it allows us to consider the problem of the finding  turning points by a new geometric point of view, namely using function $\lambda(u)$.   In particular, from this point of view the solution curve of \eqref{F} corresponds to a trough of $-\lambda(u)$ such that the continuation methods tend to tracing along it, while the geometry of the extended functional method dictates to avoid doing so.  Schematically, the geometrical difference between these approaches can be seen in  Figures \ref{fig:f2} and \ref{fig:matcont} (see \cite{matcont} for Fig. \ref{fig:matcont}).
Another advantage of the new approach is its conceptual generality and simplicity. This allows us to expect  that it can be developed in finding other types of bifurcations such as Hopf bifurcations, singularities of multidimensional parametric  and non-stationary  problems (see e.g. in \cite{BobkovIlyasov}, \cite{BobIl}, \cite{IlD} for some theoretical framework). It should be noted that to solve  maxmin problem \eqref{MiMa} (minmax problem \eqref{MaMi}) we used only one of the approaches which is based on steepest ascent direction method for piecewise smooth functions. However, there are other methods for solving such kind of problems that can be also useful (see e.g. \cite{bagir1, DemyanSPP,  kiwiel,  lemar2, fedorov, rock, shor, wolfe}). 

Although our algorithm is not yet configured to work effectively certain advantages of this approach can be seen. 
The method does not depend
on the choice of the initial point $(u_0, \lambda_0)$. 
Construction of the iterative sequence $(u^k)$ by MAQDSA consists of only one step and it does not require  additionally to implement the correction step. 
The similar systems of linear equations of the form \eqref{curav} are solved   both by MAQDSA  and by continuation methods (see Remarks \ref{DirM}, \ref{RemC}). However, by MAQDSA this system of equations (see \eqref{D}) has the dimension $N_\varepsilon(u^k)=\{i\in [1:n]:~|f_i(u^k)-\lambda(u^k)|<\varepsilon \}$ which is, in general, less then $n$. Implementation of a more detailed step-size control  of $\varepsilon$ is likely to reduce the  dimensions  of systems \eqref{D} using in the algorithm. We are currently investigating this issue to develop algorithms (based on an extended functional method) applicable to large-scale problems.

As mentioned above, the investigations presented in the current paper were not directed to the recognition   
of whether the found point $(u^*,\lambda^*)$ by  MAQDSA is really a turning point or/and a maximal turning point. However, we believe that this issue can be solved in the framework of the extended functional method and we are presently developing it in this direction. 

{\bf Acknowledgment}

Y. Il'yasov was partially supported by grant RFBR 14-01-00736-p-a. A. Ivanov was partially supported by grant RFBR 13-01-00294-p-a.

\end{document}